\documentclass[10pt]{article}

\usepackage[english]{babel}
\usepackage[utf8]{inputenc}
\usepackage{times}
\usepackage{amsmath,amsthm,amssymb,mathrsfs}
\usepackage{graphicx}
\usepackage[a4paper,bindingoffset=0.5cm,left=2cm,right=2cm,top=2cm,bottom=2cm,footskip=.8cm]{geometry}
\usepackage[md]{titlesec}
\usepackage{color,tikz}
\usepackage{url,color}
\definecolor{light}{gray}{0.8}

\newtheorem{theorem}{Theorem}[section]

\usepackage{bm}
\newcommand{\bs}{\pmb}

\newcommand{\vm}{{\bs m}}

\newtheorem{prop}{Proposition}

\theoremstyle{remark}

\newcommand{\hI}{^{({\mbox {\footnotesize {\sc i}}})}}
\newcommand{\hII}{^{({\mbox {\footnotesize {\sc ii}}})}}

\newcommand{\vb}{\vspace{3.2mm}}

\newcommand{\KK}{\kappa}
\newcommand{\cK}{\mathscr K}
\newcommand{\cL}{\mathscr L}

\title{Functional central limit theorems for\\Markov-modulated infinite-server systems}

\author{J.\ Blom$\,^\star$, K.\ De Turck$\,^\dagger$, M. Mandjes$\,^{\bullet,\star}$}

\date{\today}

\begin{document}
\maketitle

\begin{abstract} \noindent In this paper we study the Markov-modulated
    M/M/$\infty$ queue, with a focus on the correlation structure of the number
    of jobs in the system.  The main results describe the system's asymptotic
    behavior under a particular scaling of the model parameters in terms of a
    functional central limit theorem. More specifically, relying on the
    martingale central limit theorem, this result is established, covering the
    situation in which the arrival rates are sped up by a factor $N$ and the
    transition rates of the background process by $N^\alpha$, for some
    $\alpha>0$. The results reveal an interesting dichotomy, with crucially
    different behavior for $\alpha>1$ and $\alpha<1$, respectively. The
    limiting Gaussian process, which is of the Ornstein-Uhlenbeck type, is
    explicitly identified, and it is shown to be in accordance with   explicit
    results on the mean, variances and covariances of the number of jobs in the
    system.

\vspace{3mm}


\vspace{2mm}

\begin{itemize}
\item[$^\bullet$] Korteweg-de Vries Institute for Mathematics,
University of Amsterdam, Science Park 904, 1098 XH Amsterdam, the Netherlands. 
\item[$^\star$] CWI, P.O. Box 94079, 1090 GB Amsterdam, the Netherlands.
\item[$^{\dagger}$] {Laboratoire Signaux et Syst\`emes (L2S, CNRS UMR8506), \'Ecole CentraleSup\'elec, Universit\'e Paris Saclay, 3 Rue Joliot Curie, Plateau de Moulon, 91190 Gif-sur-Yvette, France.}
\end{itemize}

\noindent
M.\ Mandjes is also with  E{\sc urandom}, Eindhoven University of Technology, Eindhoven, the Netherlands, 
and IBIS, Faculty of Economics and Business, University of Amsterdam,
Amsterdam, the Netherlands. M. Mandjes' research is partly funded by the NWO Gravitation project NETWORKS, grant number 024.002.003.

\end{abstract}

\newcommand{\hN}{^{(N)}}
\newcommand{\MEAN}{{\mathbb E}}
\newcommand{\PROB}{{\mathbb P}}
\newcommand{\VAR}{{\mathbb V}{\rm ar}}
\newcommand{\COV}{{\mathbb C}{\rm ov}}

\newpage

\section{Introduction}
This paper studies the  infinite-server queue modulated by a finite-state irreducible continuous-time Markov chain $J$; 
when the so-called {\it background process} $J$ is in state $i$, jobs arrive according to a Poisson process with rate 
$\lambda_i$, while the departure rate is given by $\mu_i$.
The resulting Markov-modulated infinite-server queue has attracted some attention during the past decades; see e.g.\ the 
early contributions \cite{DAURIA, KEILSONSERVI1993, OCINNEIDEPURDUE}. In these papers the main results were in terms of 
systems of (partial) differential equations characterizing probability generating functions related to the system's 
transient behavior, and recursions enabling the evaluation of the corresponding moments. 

\vb

In a series of recent papers \cite{dave,BKMT,BKMTu,BTMASMTA,PEIS,BTM,BM}, substantial attention has been paid to the 
asymptotic behavior of Markov-modulated infinite-server queues in specific scaling regimes. 
In these parameter scalings the arrival rates are typically inflated by a factor $N$, while the transition rates of the 
background process are sped up by a factor $N^\alpha$ for some $\alpha\ge 0$. The objective is to analyze the transient 
distribution of the number of 
{jobs} in the system at time $t$, to be denoted by $M^{(N)}(t)$, in the limiting regime that $N$ grows large. 

The asymptotic results derived come in three flavors: (i)~large deviations ({\sc ld}) results{,} 
describing the tail probabilities ${\mathbb P}(M^{(N)}(t)/N\ge a)$ for $N$ large{;} 
(ii)~central-limit-theorem ({\sc clt}) type of results, describing the convergence of $M^{(N)}(t)$ 
(after centering and normalization) to a Normally distributed random variable; 
and (iii)~functional central limit theorems ({\sc fclt}\,s), describing the convergence of the process $M^{(N)}(\cdot)$ 
to an appropriate Gaussian process. 

Importantly, two model variants can be distinguished, with their own specific  departure processes. 
\begin{itemize}
\item[$\circ$]
In the first, to be referred to as Model {\sc i}, each job present is experiencing a departure rate $\mu_i$ when $J$ is
in state $i$; as a consequence, this hazard rate may change during the job's sojourn time (when the background process 
makes a transition). 
\item[$\circ$]
In the second, Model {\sc ii}, the job's sojourn time is sampled upon arrival: when the background process is then in 
state $i$, it has an exponential distribution with mean $1/\mu_i$, and hence the corresponding hazard rate is constant 
over its lifetime. 
\end{itemize}

Fig.\ \ref{figschema} summarizes the results that have been established so far. In the {\sc ld} domain, the papers
\cite{BKMTu, BTM,BM} cover, for both models,  the regime in which the background process is relatively slow  (more 
specifically, $\alpha=0$) as well as the regime in which it is essentially faster than the arrival process ($\alpha>1$). 
Also in the {\sc clt} regime the picture is complete, with results for Models {\sc i} and {\sc ii}, and with both 
slow ($\alpha<1$) and fast ($\alpha>1$) switching of the background process. In terms of {\sc fclt}\,s, however, not all 
cases are covered: the only result derived so far \cite{dave} concerns the case that $\mu_i=\mu$ for all $i$, i.e., the
case in which Models {\sc i} and {\sc ii} actually coincide; we may refer to this model as to `Model~0'. The main 
contribution of the present paper is the derivation of {\sc fclt}\,s for Models {\sc i} and {\sc ii}; this is done in 
Sections 5 and 6, respectively. These findings, with a limiting Gaussian process of 
the Ornstein-Uhlenbeck type, turn out to be in accordance with explicit expressions for means, variances, and 
covariances in these models, as we present in Sections~3~and~4. 
We conclude in Section~7 with some numerical experiments.

{\footnotesize
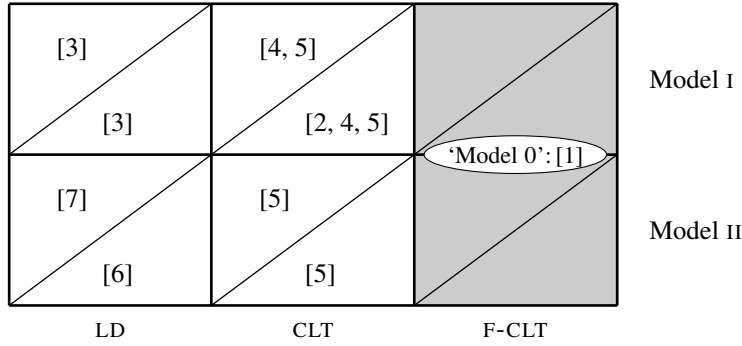
\begin{figure}
\bigskip
\begin{center}
\begin{tikzpicture}

\filldraw[fill=gray!40!white, draw=black] (7.33,-0.6) rectangle (10, 3.4);

\draw[line width=1pt] (2,-0.6)--(10,-0.6);
\draw[line width=1pt] (2,1.4)--(10,1.4);
\draw[line width=1pt] (2,3.4)--(10,3.4);
\draw[line width=1pt] (2,-0.6)--(2,3.4);
\draw[line width=1pt] (4.66,-0.6)--(4.66,3.4);
\draw[line width=1pt] (7.33,-0.6)--(7.33,3.4);
\draw[line width=1pt] (10,-0.6)--(10,3.4);

\draw[line width=0.5pt] (10,3.4)--(7.33,1.4);
\draw[line width=0.5pt] (10,1.4)--(7.33,-0.6);
\draw[line width=0.5pt] (7.33,3.4)--(4.66,1.4);
\draw[line width=0.5pt] (7.33,1.4)--(4.66,-0.6);
\draw[line width=0.5pt] (4.66,3.4)--(2,1.4);
\draw[line width=0.5pt] (4.66,1.4)--(2,-0.6);

\draw (2.5,0.8) node[right] {$\mbox{\cite{BM}}$};
\draw (2.5,2.8) node[right] {$\mbox{\cite{BKMTu}}$};
\draw (3.1,1.8) node[right] {$\mbox{\cite{BKMTu}}$};
\draw (3.1,-0.2) node[right] {$\mbox{\cite{BTM}}$};

\draw (5.166,0.8) node[right] {$\mbox{\cite{PEIS}}$};
\draw (5.166,2.8) node[right] {$\mbox{\cite{BTMASMTA,PEIS}}$};
\draw (5.766,1.8) node[right] {$\mbox{\cite{BKMT, BTMASMTA, PEIS}}$};
\draw (5.766,-0.2) node[right] {$\mbox{\cite{PEIS}}$};

\draw (10.3,2.4) node[right] {$\mbox{Model {\sc i}}$};
\draw (10.3,0.4) node[right] {$\mbox{Model {\sc ii}}$};
\draw (3.33,-0.7) node[below] {$\mbox{\sc ld}$};
\draw (6,-0.7) node[below] {$\mbox{\sc clt}$};
\draw (8.66,-0.7) node[below] {$\mbox{\sc f-clt}$};

\draw[fill=white] (8.66,1.4) ellipse (12mm and 2.5mm);
\draw (8.66,1.37)node {$\mbox{\small `Model {\sc 0}':\hspace{0.5mm}\cite{dave}}$};

\end{tikzpicture}\end{center}\caption{\label{figschema}Graphical illustration of literature on Markov-modulated 
infinite-server queues. Upper-left triangle: fast regime; lower-right triangle: slow regime. White area: regimes 
covered by earlier work; shaded areas: regimes not covered yet.}\end{figure}}

\section{Notation, preliminaries}
Let $J(t)$ denote an irreducible continuous-time Markov chain on the (finite) state space $\{1,\ldots,d\}$, with 
transition rate matrix $Q=(q_{ij})_{i,j=1}^d$ and (unique) invariant probability measure ${\boldsymbol \pi}$. 
In addition, we let $p_{ij}(t):={\mathbb P}(J(t)=j\,|\,J(0)=i).$ It is assumed that $J(0)$ is distributed according to 
${\bs{\pi}}$.

The process $J(t)$ is referred to as the {\it background process}, and regulates an infinite-server queue.
When $J(t)$  is in state $i$, jobs arrive at the queueing resource according to a Poisson process of rate $\lambda_i$. 
Regarding the way in which these jobs are handled, two variants are distinguished:
\begin{itemize}
\item[$\circ$] In Model {\sc i} the hazard rate of jobs leaving is $\mu_i$ when the background process is in state $i$. 
Observe that this hazard rate may change during the lifetime of the job, when the background process jumps. 
\item[$\circ$] In Model {\sc ii} job durations are sampled upon arrival: they are drawn from an exponential distribution
with mean $1/\mu_i$ if the background process is in state $i$ when the job enters the system.
\end{itemize}
Throughout this paper we write $\bs{\lambda}:=(\lambda_1,\ldots,\lambda_d){^{\rm T}}$ and 
${\Lambda}:={\rm diag}\{\bs{\lambda}\}$,
and likewise $\bs{\mu}:=(\mu_1,\ldots,\mu_d){^{\rm T}}$ and ${\mathcal M}:={\rm diag}\{\bs{\mu}\}$.
We also define $\lambda_\infty:={\bs \pi}^{\rm T}{\sc{\bs{\lambda}}}$ and 
$\mu_\infty:={\bs \pi}^{\rm T}{\sc{\bs{\mu}}}$. 

In Sections 3 and 4 we consider explicit expressions for the means, variances and covariances in the 
{\it unscaled system}. There we denote by $M(t)$ the number of jobs present at time $t$, for $t\ge 0.$ For simplicity, it is assumed that the system starts empty at time $0$, i.e., $M(0)=0$. 

In Sections 3 and 4 we also analyze the obtained expressions for the mean, variance and covariance in a specific parameter scaling, viz.\ we replace the arrival rates $\bs{\lambda}$ by $N\bs{\lambda}$, and the generator matrix $Q$ by $N^\alpha Q$, for some $\alpha>0$, and let $N$ grow large. It is in this asymptotic regime that we also establish our {\sc fclt}\,s in Sections 5 and 6. For these scaled models we write $M^{(N)}(t)$ for the number of jobs at time $t$, to emphasize the dependence on the scaling parameter. 

\vb

In the sequel, we use the concept of {\it deviation matrices.} Define the $(i,j)$-th element of the {\it exponentially weighted deviation matrix} $D^{({\boldsymbol \gamma})}$, as a function of the vector ${\bs\gamma}\in{\mathbb R}_+^d$, by \[D^{({\boldsymbol \gamma})}_{ij} := \int_0^\infty e^{-\gamma_i t} \left(p_{ij}(t)-\pi_j\right){\rm d} t.\] The matrix $D:= D^{({\boldsymbol 0})}$ is the canonical {\it deviation matrix}.  In the sequel, also the matrix $\Pi:={\bs 1}{\bs \pi}^{\rm T}$ plays a role, as well as the fundamental matrix $F:=D+\Pi.$ A number of identities hold: $QF=FQ=\Pi-I$, $\Pi F=F\Pi=\Pi$, and $F{\bs 1}=\Pi{\bs 1}={\bs 1}.$

\section{Mean, variance, and covariance for model {\sc i}}

The first part of this section presents explicit formulae for the mean, variance, and covariance in Model {\sc i}. In the second part these turn out to allow for a more explicit characterization in particular asymptotic regimes.

\subsection{Explicit formulae}

Our goal is to devise a method to compute ${\mathbb C}{\rm ov}(M(t),M(t+u))$. To this end, 
the object that we study first is, for $u\ge 0$ fixed, the bivariate probability generating function
\[ \Xi_{ij}(z,w,t,u):=\MEAN\left( z^{M(t)}w^{M(t+u)}\, 1{\{J(t)=i,J(t+u)=j\}}\right),\]
which implicitly contains all information about the joint distribution of $M(t)$ and $M(t+u)$.
In matrix notation, we obtain in Appendix~A, suppressing the arguments for ease of notation,
\begin{equation}
\label{xi}\frac{\partial\Xi }{\partial t}  = (z-1) \Lambda\Xi+(w-1)\Xi\Lambda
-(z-1){\mathcal M}\frac{\partial \Xi}{\partial z}  - 
(w-1)\frac{\partial \Xi}{\partial w} {\mathcal M}+Q^{\rm T}\Xi+\Xi Q.
\end{equation}
We now point out how to compute the covariance between $M(t)$ and $M(t+u)$ from this system of partial differential equations.

\vb

To this end, we first define the three matrices
\begin{align}
  & E(t,u)\equiv(E_{ij}(t,u))_{i,j=1}^d,\:\:\:\:\text{ where } E_{ij}(t,u) := {\mathbb E} M(t) 1\{J(t)=i,J(t+u)=j\} \nonumber\\
  & G(t,u)\equiv(G_{ij}(t,u))_{i,j=1}^d,\:\:\:\:\text{ where } G_{ij}(t,u) := {\mathbb E} M(t+u) 1\{J(t)=i,J(t+u)=j\} \nonumber\\
  & C(t,u)\equiv(C_{ij}(t,u))_{i,j=1}^d,\:\:\:\:\text{ where } C_{ij}(t,u) := {\mathbb E} M(t) M(t+u) 1\{J(t)=i,J(t+u)=j\} \nonumber\\
\end{align}

It follows from the moment-generating property of generating functions that 
\begin{equation}E_{ij}(t,u) =\lim_{z,w\uparrow 1} \frac{\partial \Xi_{ij} }{\partial z} ,\:\:\:\: 
G_{ij}(t,u) =\lim_{z,w\uparrow 1} \frac{\partial \Xi_{ij}}{\partial w},\:\:\:\: 
C_{ij}(t,u) =\lim_{z,w\uparrow 1} \frac{\partial^2 \Xi_{ij}}{\partial z\partial w}.
\end{equation}
From the partial differential equation (\ref{xi}) that defines $\Xi$, we can find the following systems of ordinary differential equations for the matrices $E(t,u)$, $G(t,u)$ and $C(t,u)$.
We demonstrate how this is done for the equation involving $E(t,u).$ 
Differentiate (\ref{xi}) with respect to $z$, and take the limit of $w,z\uparrow 1$. 
Recalling that $J(0)$ is distributed according to ${\bs \pi}$, it is straightforward to obtain, with $K_{ij}(u):=\pi_i p_{ij}(u)$,
\[E'(t,u) = \Lambda K(u) -{\mathcal M}E(t,u) +Q^{\rm T} E(t,u)+E(t,u)Q,\]
where the derivative in the left-hand side is again with respect to $t$.
We can derive the {\sc ode}s for $G(t,u)$ in the same manner,
    \[G'(t,u)=K(u)\Lambda-G(t,u){\mathcal M}+Q^{\rm T} G(t,u)+G(t,u)Q.\]
Similarly, for $C(t,u)$ we have:
    \[C'(t,u) =\Lambda G(t,u)+E(t,u)\Lambda -{\mathcal M}C(t,u)-C(t,u){\mathcal M}+Q^{\rm T}C(t,u)+C(t,u)Q,\]

\vb

The above differential equations are {\it matrix-valued} systems of linear differential equations, which can be converted into {\it  vector-valued} systems of linear differential equations,  relying on the concept of `vectorization'. We show this idea for the matrix $E(t,u).$ We take the columns of $E(t,u)$, and put them into a vector ${\bs e}(t,u)$ of dimension $d^2$, such that the first $d$ entries are $E_{11}(t,u)$ up to $E_{d1}(t,u)$, entries $d+1$ up to $2d$ correspond to $E_{12}(t,u)$ up to $E_{d2}(t,u)$, etc.; we write ${\bs e}(t,u):={\rm vec} (E(t,u))$. Likewise, ${\bs g}(t,u):={\rm vec} (G(t,u))$, ${\bs c}(t,u):={\rm vec} (C(t,u))$ and ${\bs k}(u):={\rm vec}(K(u))$.

For $d\times d$ matrices $A$, $B$, and $C$, and with {as usual} $A\otimes B$  denoting the Kronecker 
product and {$A\oplus B:= A \otimes I + I \otimes B$ the Kronecker sum} of the matrices $A$ and $B$, recall
\[{{\rm vec}(AB) = (I\otimes A){\rm vec}(B) = (B^{\rm T}\otimes I){\rm vec}(A)}
\:\:\:\mbox{and}\:\:\:\:
{\rm vec}(ABC) = (C^{\rm T}\otimes A){\rm vec}(B),\]
with $I$ the $d\times d$ identity matrix. 
We thus obtain the following equations in terms of Kronecker sums and products:
\[{\bs e}'(t,u)= 
{(I\otimes \Lambda)}\bs{k}(u)-(I\otimes {\mathcal M}){\bs e}(t,u)+  (Q^{\rm T}\oplus Q^{\rm T}){\bs e}(t,u).\]
An equation for ${\bs g}(t,u)$ can be found analogously:
\[{\bs g}'(t,u)= {( \Lambda\otimes I)}\bs{k}(u)
-( {\mathcal M}\otimes I){\bs g}(t,u)+ (Q^{\rm T}\oplus Q^{\rm T}){\bs g}(t,u).\]
 Along the same lines we obtain
\[{\bs c}'(t,u) = {(I\otimes \Lambda)} {\bs{g}}(t,u)+{(\Lambda \otimes I)} {\bs{e}}(t,u)
 -({\mathcal M}\oplus{\mathcal M}) {\bs c}(t,u) +  (Q^{\rm T}\oplus Q^{\rm T}){\bs c}(t,u),\]
the derivatives in the left-hand sides being again with respect to $t$.
Observe that $Q\oplus Q$ is again a transition rate matrix, and ${\mathcal M}\oplus{\mathcal M}$ a diagonal matrix with non-negative entries.

The systems describing ${\bs e}(t,u)$ {and} ${\bs g}(t,u)$ are  standard systems of non-homogeneous linear differential equations, which can be solved with standard techniques. Then the solution can be plugged into the differential equation describing ${\bs c}(t,u)$, which is then also a system of non-homogeneous linear differential equations. We summarize the results in the following proposition.

\begin{prop}
    The matrix-valued functions $E(t,u), G(t,u)$ and $C(t,u)$ satisfy the following {\sc ode}s:
    \[E'(t,u) = \Lambda K(u) -{\mathcal M}E(t,u) +Q^{\rm T} E(t,u)+E(t,u)Q,\]
    \[G'(t,u)=K(u)\Lambda-G(t,u){\mathcal M}+Q^{\rm T} G(t,u)+G(t,u)Q.\]
    \[C'(t,u) =\Lambda G(t,u)+E(t,u)\Lambda -{\mathcal M}C(t,u)-C(t,u){\mathcal M}+Q^{\rm T}C(t,u)+C(t,u)Q,\]

Moreover, the vectorized versions ${\bs e}(t,u)$, ${\bs g}(t,u)$ and ${\bs c}(t,u)$, of the matrices $E(t,u)$, $G(t,u)$ and $C(t,u)$ satisfy the following linear differential equations.
\[{\bs e}'(t,u)= {(I\otimes \Lambda)}\bs{k}(u)-(I\otimes {\mathcal M}){\bs e}(t,u)+  (Q^{\rm T}\oplus Q^{\rm T}){\bs e}(t,u).\]
\[{\bs g}'(t,u)= {( \Lambda\otimes I)}\bs{k}(u)
-( {\mathcal M}\otimes I){\bs g}(t,u)+ (Q^{\rm T}\oplus Q^{\rm T}){\bs g}(t,u).\]
\[{\bs c}'(t,u) = {(I\otimes \Lambda)} {\bs{g}}(t,u)+{(\Lambda \otimes I)} {\bs{e}}(t,u)
 -({\mathcal M}\oplus{\mathcal M}) {\bs c}(t,u) +  (Q^{\rm T}\oplus Q^{\rm T}){\bs c}(t,u).\]

 All occurring derivatives are with respect to $t$.
\end{prop}

\vb

We have now devised a procedure to compute the covariance ${\mathbb C}{\rm ov}(M(t),M(t+u))$.  To this end, first realize that, with $e(t):= {\mathbb E} M(t)$ and ${\bs 1}$ denoting here a $d^2$-dimensional all-ones vector,
\[e(t)={\bs 1}^{\rm T} {\bs e}(t,u),\:\:\:
e(t+u)= {\bs 1}^{\rm T} {\bs g}(t,u).\]
As a consequence,
\[{\mathbb C}{\rm ov}(M(t),M(t+u)) = {\bs 1}^{\rm T} {\bs c}(t,u) -{e}(t)\,{e}(t+u).\]

\subsection{Two specific limiting regimes}\label{SR}

In this subsection, we consider two particular limiting regim{e}s, in which the expressions simplify considerably.

\vb

$\rhd$\:
Let us first consider the behavior for $t\to\infty$. It is readily verified that
\[{\bs e}(\infty,u){:=}\lim_{t\to\infty}{\bs e}(t,u) = \left( (I\otimes{\mathcal M}) -
(Q^{\rm T}\oplus Q^{\rm T})  \right)^{-1}{(I\otimes \Lambda)}\bs{k}(u) ,\]
and hence
\[e(\infty)= \bs{1}^{\rm T}{\bs e}(\infty,u) =\bs{1}^{\rm T} \left( (I\otimes{\mathcal M}) -(Q^{\rm T}\oplus Q^{\rm T})  
\right)^{-1}{(I\otimes \Lambda)}\bs{k}(u) 
= \bs{1}^{\rm T}{\bs g}(\infty,u) .\]
For $u=0$ we obtain the solution from O'Cinneide and Purdue \cite[Thm.\ 3.1]{OCINNEIDEPURDUE}. 

Along the same lines,
\[{\bs c}(\infty,u){:}=\lim_{t\to\infty}{\bs c}(t,u)=\left( ({\mathcal M}\oplus{\mathcal M})  -  
(Q^{\rm T}\oplus Q^{\rm T})\right)^{-1} ({(I\otimes \Lambda)}{\bs{g}}(\infty,u)+
{(\Lambda \otimes I)}{\bs e}(\infty,u)) .\]
We have thus derived an expression for the limit of ${\mathbb C}{\rm ov}(M(t),M(t+u))$ as $t\to\infty$:
\[\lim_{t\to\infty}{\mathbb C}{\rm ov}(M(t),M(t+u)) = {\bs 1}^{\rm T} {\bs c}(\infty,u) - (e(\infty))^2.\]

\vb

$\rhd$\, Next, we consider the following scaling: we replace ${\boldsymbol\lambda} \mapsto N {\boldsymbol\lambda}$ and $Q\mapsto N^\alpha Q$, for $\alpha>0$.  In this regime, the pace with which the arrival process is sped up, differs from that corresponding to the background process. As we will see below, the  situation $\alpha>1$ crucially differs from $\alpha<1$; this was already observed earlier in e.g.\ {\cite{dave,PEIS}}.  As mentioned before, to stress the dependence on $N$, we write $M\hN(t)$ rather than $M(t).$ It is this scaling that is imposed in Section 5, and under which an {\sc fclt} is established. We now identify the associated mean and (co-)variance, relying on elementary techniques.

Let ${\bs m}\hN(t)\equiv {\bs m}(t)$ the $d$-dimensional {\it row-}vector, with ${\mathbb E} M\hN(t)1{\{J(t)=i\}}$ on the $i$-th position.  According to \cite[Thm.\ 3.2]{OCINNEIDEPURDUE}, ${\bs m}(t)$ satisfies the following non-homogeneous linear differential equation:
\begin{equation}\label{mteq}{\bs\pi}^{\rm T}N\Lambda - {\bs m}(t) ({\mathcal M}-N^\alpha Q) = {\bs m}'(t).\end{equation}
In  {\cite{PEIS}} we proved that, with $\varrho\hI:=\lambda_\infty/\mu_\infty$, 
\begin{equation}\label{lime}{\mathbb E} M\hN(t) = N\,\varrho\hI\,(1-e^{-\mu_\infty t})  +o(N).\end{equation}
Now define $\varrho\hI(t):= \varrho\hI(1-e^{-\mu_\infty t})$ and
\begin{equation}\label{VS}\varsigma\hI(t):=2\int_0^t e^{-2\mu_\infty (t-s)} {\bs\pi}^{\rm T}\left(\Lambda-
{\mathcal M}\varrho\hI(s)\right)D\left(\Lambda-{\mathcal M}\varrho\hI(s)\right){\bs 1}\,{\rm d}s .
\end{equation}
In Appendix B it is shown that
\begin{equation}
\label{limc}\lim_{N\to\infty}\frac{\COV(M\hN(t),M\hN(t+u))}{N^{\max\{1,2-\alpha\}}}= v\hI(t,u):=e^{-\mu_\infty u} 
\left(\varsigma\hI(t)1_{\{\alpha\le 1\}} +\varrho\hI(t)1_{\{\alpha\ge 1\}}\right).\end{equation}
We conclude that under this parameter scaling the covariance exhibits the same dichotomy as the one 
observed in  {\cite{PEIS}} for the variance, i.e., behaving crucially different for $\alpha<1$ and $\alpha>1$.
In the latter regime, the system essentially behaves as a (non-modulated) M/M/$\infty$ queue, with arrival rate $\lambda_\infty$ and service rate $\mu_\infty$, whereas for $\alpha<1$ the full transition rate matrix $Q$ plays a role (as $\varsigma\hI(t)$ involves the deviation matrix $D$).

\section{Mean, variance, and covariance for model {\sc ii}}

As we saw above, for Model {\sc i} the mean, variance and covariance can be determined by solving specific non-homogeneous linear differential equations; for Model {\sc ii}, however, the analysis is simpler, and can be performed by relying on the law of total (co-)variance, as shown in Section \ref{EF}. Focusing on the same limiting regimes as we have studied for Model {\sc i}, the expressions become more explicit; see Section \ref{TSLR}. 

\subsection{Explicit formulae}\label{EF}
The mean of $M(t)$ for Model {\sc ii} was already determined in \cite{BKMT}; recalling from e.g.\ \cite{DAURIA} the observation that $M({t})$ obeys a Poisson distribution with the random parameter 
$\MEAN(M({t})\,|\,J)$, we conclude that
\[{\mathbb E}M(t) = \MEAN\left(\MEAN(M(t)\,|\,J) \right) = {\mathbb E}\left(
\int_0^t \lambda_{J(s )} e^{-\mu_{J( s)}(t-s)}{\rm d}s
\right)=\sum_{i=1}^d \pi_i \frac{\lambda_i}{\mu_i}\left(1-e^{-\mu_i t}\right)=:\varrho\hII(t),\]
with $J\equiv (J( s))_{s=0}^t$.

Now concentrate on the evaluation of the covariance between $M(t)$ and $M(t+u)$; 
assume, without loss of generality, that $u\ge 0$.
The `law of total covariance' entails that
\begin{equation}\label{TC}
\COV(M(t),M(t+u)) = \MEAN(\COV(M(t),M(t+u)\,|\,J))+\COV(\MEAN(M(t)\,|\,J),\MEAN(M(t+u)\,|\,J)).
\end{equation}
In Appendix D, we evaluate both terms, so as to obtain
\begin{equation}\label{eqcov}
\COV(M(t),M(t+u))=\sum_{i=1}^d \pi_i \frac{\lambda_i}{\mu_i}\left(1-e^{-\mu_i t}\right)e^{-\mu_i u}
+{\bs \lambda}^{\rm T}\cK(t,u){\bs \lambda}+{\bs \lambda}^{\rm T}\cL(t,u){\bs \lambda};
\end{equation}
the precise form of the matrices $\cK(t,u)$ and $\cL(t,u)$ is given in Appendix D as well.

\subsection{Two specific limiting regimes}\label{TSLR}
In this subsection, we consider the two particular limiting regimes that we studied earlier, in Section~\ref{SR}, for Model~{\sc i}. As it turns out, in these regimes the expressions simplify considerably.

\vb

$\rhd$\, In the first regime, we consider $\COV(M(t),M(t+u))$ for $t\to\infty$. Going through the calculations, relying 
on the explicit expressions for $\cK(t,u)$ and $\cL(t,u)$ as given in Appendix D, we obtain
\begin{eqnarray*}
\lim_{t\to\infty}
\COV(M(t),M(t+u))&=&\sum_{i=1}^d \pi_i \frac{\lambda_i}{\mu_i}e^{-\mu_iu}+
\sum_{i=1}^d\sum_{j=1}^d \pi_i\frac{\lambda_i\lambda_j}{\mu_i+\mu_j} e^{-\mu_j u} D^{({\boldsymbol\mu})}_{ij}\\
&&
+\,\sum_{i=1}^d\sum_{j=1}^d \pi_j
\frac{\lambda_i\lambda_j}{\mu_i+\mu_j} \int_0^u e^{-\mu_j u+\mu_i w} (p_{ji}(w)-\pi_i){\rm d}w\\
&&+\,
\sum_{i=1}^d\sum_{j=1}^d \pi_j
\frac{\lambda_i\lambda_j}{\mu_i+\mu_j} \int_u^\infty e^{\mu_i u-\mu_j w} (p_{ji}(w)-\pi_i){\rm d}w,
\end{eqnarray*}
also entailing that
\[\lim_{t\to\infty}
\VAR M(t) = \sum_{i=1}^d \pi_i \frac{\lambda_i}{\mu_i} + 2 
\sum_{i=1}^d\sum_{j=1}^d \pi_i\frac{\lambda_i\lambda_j}{\mu_i+\mu_j} D^{({\boldsymbol\mu})}_{ij}.\]

\vb

$\rhd$\, In the second limit, we replace ${\boldsymbol\lambda} \mapsto N {\boldsymbol\lambda}$ and 
$Q\mapsto N^\alpha Q$, for $\alpha>0.$ The {\sc fclt} under this scaling is proven in Section 6; we here find  the 
corresponding mean and (co-)variance.
It turns out that, for $N$ large, \[
\COV\left(M\hN(t),M\hN(t+u)\right)\sim N\sum_{i=1}^d e^{-\mu_iu}\varrho_i\hII(t)+
N^{2-\alpha}\sum_{i=1}^de^{-\mu_iu}\varsigma_i\hII(t),\]
with $\varrho_i\hII  :=
\pi_i\lambda_i/\mu_i$ and $\varrho_i\hII(t):=\varrho_i\hII  \cdot(1-e^{-\mu_i t})$
and
\[\varsigma_i\hII(t):=\sum_{j=1}^d \frac{\lambda_i\lambda_j}{\mu_i+\mu_j}  
 \left(1-e ^{-(\mu_i+\mu_j) t}\right) (\pi_jD_{ji}+ \pi_iD_{ij} ).\]
 We conclude that
 \begin{equation}
\label{limcII}
\lim_{N\to\infty}\frac{\COV(M\hN(t),M\hN(t+u))}{N^{\max\{1,2-\alpha\}}}= v\hII(t,u):=
\sum_{i=1}^de^{-\mu_i u} \left(\varsigma_i\hII(t)1_{\{\alpha\le 1\}} +\varrho_i\hII(t)1_{\{\alpha\ge 1\}}\right).
\end{equation}
 We observe that the same dichotomy applies as the one we have observed for Model {\sc i}:
for $\alpha>1$ the number of jobs in the system behaves `Poissonian', with mean and variance scaling essentially linearly
with $N$, both with proportionality constant $\varrho\hII(t)$. For $\alpha<1$, as seen earlier in 
e.g.\ {\cite{PEIS}}, the 
variance grows superlinearly with $N$, with a proportionality constant that involves the deviation matrix $D$.

\section{Functional central limit theorem for Model {\sc i}}
In Section \ref{SR} we considered the covariance of the number of jobs in the system under a specific 
scaling: ${\boldsymbol\lambda} \mapsto N {\boldsymbol\lambda}$ and $Q\mapsto N^\alpha Q$, for $\alpha>0$. 
In this section, we prove that for a given $t$ the random variable $M^{(N)}(t)$ obeys a central limit theorem; moreover, 
we prove the stronger property that after centering and normalizing the process $M^{(N)}(t)$, there is weak convergence 
to a specific Gaussian process. We essentially adopt the methodology used in \cite{mmou}; some steps that are fully 
analogous to those in \cite{mmou} are described concisely. In the sequel, we let $Z^{(N)}_i(t)$ be the indicator function 
of the event $\{J^{(N)}(t)=i\}$, where $J^{(N)}(t)$ is a Markov chain with transition rate matrix $N^\alpha Q$.

First observe that, with $P_1(\cdot)$ and $P_2(\cdot)$ two independent unit-rate Poisson processes, it is straightforward 
to see that $M^{(N)}(t)$ can be written as
\begin{equation}\label{diffpp}
M^{(N)}(t)= P_1\left(N \int_0^t \sum_{i=1}^d \lambda_i Z_i^{(N)}(s){\rm d}s\right)- P_2
\left(\int_0^t \sum_{i=1}^d \mu_i M^{(N)}(s) Z_i^{(N)}(s){\rm d}s\right).
\end{equation}
Now impose the following centering and normalization, with $\beta:= \max\{1, 2-\alpha\}{/2}$,
\[\tilde M^{(N)}(t) := N^{-\beta} \left(M^{(N)}(t) -N \varrho\hI(t)\right),\]
where $\varrho\hI(t):= \varrho\hI (1-e^{-\mu_\infty t});$
the objective of this section is to establish the convergence of $\tilde M^{(N)}(\cdot) $ to a specific Gaussian process,
essentially relying on the martingale central limit theorem; see for background on the martingale central limit theorem
e.g.\   \cite{JS, WHITT}. 

It is first realized that, as a direct implication of (\ref{diffpp}), for some martingale $\KK^{(N)}(\cdot)$,\[{\rm d} M^{(N)}(t)= N {\boldsymbol\lambda}^{\rm T} {\boldsymbol Z}^{(N)}(t){\rm d}t - 
{\boldsymbol\mu}^{\rm T} {\boldsymbol Z}^{(N)}(t)\, M^{(N)}(t){\rm d}t +{\rm d} \KK^{(N)}(t).\]
Then we rewrite this equation in terms of one for $\tilde M^{(N)}(t)$:
\begin{eqnarray*}{\rm d} \tilde M^{(N)}(t)&=& N^{1-\beta} {\boldsymbol\lambda}^{\rm T} {\boldsymbol Z}^{(N)}(t){\rm d}t - 
N^{-\beta}{\boldsymbol\mu}^{\rm T} {\boldsymbol Z}^{(N)}(t)\, M^{(N)}(t){\rm d}t +N^{-\beta}{\rm d} \KK^{(N)}(t)
-N^{1-\beta}\left(\varrho\hI\right)'\hspace{-0.5mm}(t){{\rm d}t}\\
&=&N^{1-\beta} {\boldsymbol\lambda}^{\rm T} {\boldsymbol Z}^{(N)}(t){\rm d}t - 
{\boldsymbol\mu}^{\rm T} {\boldsymbol Z}^{(N)}(t)\, \tilde M^{(N)}(t){\rm d}t 
-
N^{1-\beta}{\boldsymbol\mu}^{\rm T} {\boldsymbol Z}^{(N)}(t) \varrho\hI(t){{\rm d}t}\\
&&\hspace{0.02cm}
+\,N^{-\beta}{\rm d} \KK^{(N)}(t)
-N^{1-\beta}\left(\varrho\hI\right)'\hspace{-0.5mm}(t){{\rm d}t}.
\end{eqnarray*}
Following the ideas of \cite{mmou}, we now introduce
\[Y^{(N)}(t):= \exp\left({\boldsymbol\mu}^{\rm T} {\boldsymbol \zeta}^{(N)}(t)\right) \tilde{M}^{(N)}(t),\:\:\:
\mbox{where}\:\:\:\:
{\boldsymbol \zeta}^{(N)}(t):=\int_0^t {\boldsymbol Z}^{(N)}(s){\rm d}s.\]
It thus follows that, using standard stochastic differentiation rules,
\[{\rm d}Y^{(N)}(t)=  \exp\left({\boldsymbol\mu}^{\rm T} {\boldsymbol \zeta}^{(N)}(t)\right)\hspace{-1mm}
\left( N^{1-\beta}\left( {\boldsymbol\lambda}-{\boldsymbol\mu}\varrho\hI(t)\right)^{\rm T} 
{\boldsymbol Z}^{(N)}(t){\rm d}t 
+N^{-\beta}{\rm d} \KK^{(N)}(t)
-N^{1-\beta}\left(\varrho\hI\right)'\hspace{-0.5mm}(t){{\rm d}t}
\right)\hspace{-0.8mm}.\]
Now observe that, from the definition of the function $\varrho\hI(t)$, we find that 
\[\left( {\boldsymbol\lambda}-{\boldsymbol\mu}\varrho\hI(t)\right)^{\rm T}{\boldsymbol \pi}=\lambda_\infty
e^{-\mu_\infty t}=\left(\varrho\hI\right)'\hspace{-0.5mm}(t),\] 
and hence it is obtained that
\[ {\rm d}Y^{(N)}(t)=\exp\left({\boldsymbol\mu}^{\rm T} {\boldsymbol \zeta}^{(N)}(t)\right)
\left( N^{1-\beta}\left( {\boldsymbol\lambda}-{\boldsymbol\mu}\varrho\hI(t)\right)^{\rm T}
\left( {\boldsymbol Z}^{(N)}(t)-{\boldsymbol \pi}\right){\rm d}t 
+N^{-\beta}{\rm d} \KK^{(N)}(t)
\right).\]
We now analyze the two terms in the previous display separately. 
\begin{itemize}
\item[$\circ$]We first concentrate on the first term.
In \cite{mmou}, relying on the methodology developed in \cite{JS}, it was shown that the following weak convergence holds:
\[\int_0^\cdot N^{\alpha/2} \exp\left({\boldsymbol\mu}^{\rm T} {\boldsymbol \zeta}^{(N)}(s)\right)
\left( {\boldsymbol\lambda}-{\boldsymbol\mu}\varrho\hI(s)\right)^{\rm T}
\left( {\boldsymbol Z}^{(N)}(s)-{\boldsymbol \pi}\right){\rm d}s\to\int_0^\cdot e^{\mu_\infty s} {\rm d}\mathcal{G}(s),\]
where the stochastic process $\mathcal{G}(\cdot)$ is such that 
\begin{eqnarray*}\langle \mathcal{G}\rangle_t = {V}(t) &:=& 
\int_0^t \left({\boldsymbol\lambda}-{\boldsymbol\mu}\varrho\hI(s)\right)^{\rm T}
\left({\rm diag}\{{\boldsymbol \pi}\} D+D^{\rm T} {\rm diag}\{{\boldsymbol \pi}\}\right)
\left({\boldsymbol\lambda}-{\boldsymbol\mu}\varrho\hI(s)\right){\rm d}s\\
&=&2\int_0^t {\boldsymbol \pi}^{\rm T}\left({\Lambda}-{\mathcal M}\varrho\hI(s)\right)D
\left({\Lambda}-{\mathcal M}\varrho\hI(s)\right) {\boldsymbol 1}\,{\rm d}s
;\end{eqnarray*}
cf.\ Eqn.\ (\ref{VS}).
(It is noted that in \cite{mmou} the background process was sped up by a factor $N$ rather than $N^\alpha$; this 
explains that there the growth rate $\sqrt{N}$ was found, while in our setup we have $N^{\alpha/2} $.) 

Importantly, from the above we conclude that the full first term in ${\rm d}Y^{(N)}(t)$ behaves essentially proportional 
to $N^{1-\beta-\alpha/2},$ which  converges to a constant if $\alpha\leq 1$, and vanishes otherwise.
\item[$\circ$] We now consider the second term.
We note that, recalling the fact that $P_1(\cdot)$ and $P_2(\cdot)$ are independent unit-rate Poisson processes in combination with standard properties for pure jump processes,
\[\frac{\rm d}{{\rm d}t} \langle \KK^{(N)}\rangle_t = 
N {\boldsymbol\lambda}^{\rm T} {\boldsymbol Z}^{(N)}(t) + 
{\boldsymbol\mu}^{\rm T} {\boldsymbol Z}^{(N)}(t)\, M^{(N)}(t),\]
and consequently 
\[\frac{1}{N} \langle \KK^{(N)}\rangle_t = 
\int_0^t {\boldsymbol\lambda}^{\rm T} {\boldsymbol Z}^{(N)}(s){\rm d}s +
\int_0^t {\boldsymbol\mu}^{\rm T} {\boldsymbol Z}^{(N)}(s)\frac{ M^{(N)}(s)}{N}{\rm d}s.\]
Using the ergodic theorem, the first integral in the right-hand side of the previous display converges to 
${\boldsymbol\lambda}^{\rm T} {\boldsymbol \pi}\cdot t=\lambda_\infty t$. Likewise, the second integral converges to
\[\lim_{N\to\infty} \frac{1}{N} \int_0^t \sum_{i=1}^d \mu_i \,{\mathbb E} \left(M^{(N)}(s)1\{J(s)=i\}\right){\rm d}s,\]
which, due to arguments similar to those underlying (\ref{lime}), turns out to equal
\[\int_0^t \sum_{i=1}^d \mu_i \pi_i \varrho\hI(1-e^{-\mu_\infty s}){\rm d}s.\]
Hence $N^{-1}\langle \KK^{(N)}\rangle_t$ converges, as $N\to\infty$,  {to}
\[W(t):=\lambda_\infty t + \int_0^t \mu_\infty \varrho\hI(1-e^{-\mu_\infty s}){\rm d}s.\]
We conclude from the above that $\KK^{(N)}(\cdot)/\sqrt{N}$ converges to an appropriately scaled Brownian motion. 

In addition, this second term in ${\rm d}Y^{(N)}(t)$ is essentially proportional to $N^{1/2-\beta}$, 
i.e., converging to a constant if $\alpha\geq 1$, and vanishes otherwise.
\end{itemize}

{Summarizing, we have that $Y^{(N)}(t)$ converges weakly to a process $Y(t)$
which is the solution to the following stochastic differential equation:
\[{\rm d} Y(t) = {\sqrt{V'(t)1_{\{\alpha\le 1\}}+W'(t)1_{\{\alpha\ge 1\}}}}\,{\rm d}B(t),\]
where we used the property that for a standard Brownian motion $B$ and
a differentiable function $f$, we have that $B(f(t))$ is equal in distribution to
$\sqrt{f'(t)} \hat B(t)$, where $\hat B$ denotes another Brownian motion, but with the same distribution. Also, due to the ergodic theorem we have that
$\exp\left({\boldsymbol\mu}^{\rm T} {\boldsymbol \zeta}^{(N)}(t)\right)$ converges to $\exp(\mu_\infty t)$. From the definition of $Y^{(N)}(t)$,} we thus conclude the following weak convergence:  $\tilde M^{(N)}(\cdot)\to \tilde M(\cdot),$
where $\tilde M(\cdot)$ solves the stochastic differential equation
\[{\rm d}\tilde M(t) =-\mu_\infty \tilde M(t) {\rm d}t +
 {\sqrt{V'(t)1_{\{\alpha\le 1\}}+W'(t)1_{\{\alpha\ge 1\}}}}\,{\rm d}B(t),\]
for a standard Brownian motion $B(\cdot).$ Its solution is that the limiting process $\tilde M(\cdot)$ is a centered 
Gaussian process of the Ornstein-Uhlenbeck type, characterized by its covariance $v\hI(t,u)$, as given in (\ref{limc}).
\begin{theorem} Consider Model {\sc i}.  As $N\to\infty$, the process $\tilde M^{(N)}(\cdot)$ converges weakly to a 
centered Gaussian process, with covariance structure $v\hI(\cdot,\cdot)$
given in \eqref{limc}.
\end{theorem}

\section{Functional central limit theorem for Model {\sc ii}}We now shift our attention from Model {\sc i} to 
Model {\sc ii}. Essentially the same approach can be followed, with an 
important difference being that now one has to keep track of the number of jobs present {\it of each type}, to be 
denoted 
by $M^{(N)}_i(t)$ for type $i$, where `type' refers to the state the background process was in upon arrival of the job. 
We use an approach similar to the one used in the previous section, but it is noted that a viable alternative is to adapt
the approach followed in \cite{dave} for the case that the departure rates are state-independent, to that of 
Model {\sc ii}.

As in the previous section, we start by writing the $M^{(N)}_i(t)$, for $i=1,\ldots,d$ in terms of unit-rate Poisson 
processes; in self-evident notation, we now have
\[M^{(N)}_i(t) = P_{1,i}\left(N \int_0^t \lambda_i Z_i^{(N)}(s){\rm d}s\right)-
P_{2,i}\left(\int_0^t \mu_i M_i^{(N)}(s){\rm d}s\right).\]
As before, we apply centering and normalization, in that we will study, 
recalling that $\varrho_i\hII  :=
\pi_i\lambda_i/\mu_i$ and $\varrho_i\hII(t):=\varrho_i\hII  \cdot(1-e^{-\mu_it})$, 
\[\tilde M^{(N)}_i(t):= N^{-\beta}\left(M_i^{(N)}(t)-N \varrho\hII_i(t)\right),\]
where, as in the previous section, $\beta:= \max\{1, 2-\alpha\}{/2}$.
Also we have that, for martingales $\KK_i^{(N)}(t)$ (with $i=1,\ldots,d$),
\[{\rm d}M^{(N)}_i(t) =N\lambda_iZ_i^{(N)}(t) {\rm d}t - \mu_iM_i^{(N)}(t) {\rm d}t+{\rm d}\KK_i^{(N)}(t),\]
which we can express in terms of ${\rm d}\tilde M^{(N)}_i(t)$:
\begin{eqnarray*}
{\rm d}\tilde M^{(N)}_i(t) &=&N^{-\beta} {\rm d}M_i^{(N)}-N^{1-\beta}
\left(\varrho_i\hII\right)'\hspace{-1mm}(t)\,{\rm d}t\\
&=&N^{1-\beta} \lambda_i Z_i^{(N)}(t){\rm d}t -\mu_i \tilde M_i^{(N)}(t) {\rm d}t
-N^{1-\beta}\mu_i \varrho\hII_i(t){\rm d}t\\
&&+\, N^{-\beta}{\rm d}\KK_i^{(N)}(t) -N^{1-\beta}\left(\varrho_i\hII\right)'\hspace{-1mm}(t)\,{\rm d}t.
\end{eqnarray*}
Using the definition of $\varrho\hII_i(t)$, after some calculus we eventually obtain the stochastic differential
equation
\[{\rm d}\tilde M^{(N)}_i(t)=-\mu_i \tilde M_i^{(N)}(t) {\rm d}t+ N^{1-\beta} \lambda_i
\left(Z_i^{(N)}(t)-\pi_i\right){\rm d}t +N^{-\beta} {\rm d}\KK_i^{(N)}(t) .\]
Mimicking the ideas used in the previous section, we study the last two terms appearing in the right-hand side of the 
previous display separately.
\begin{itemize}
\item[$\circ$]We first concentrate on the middle term of the right-hand side of the previous display. To this end, we 
define
\[
I_i^{(N)}(t):=\int_0^t\left(Z_i^{(N)}(s)-\pi_i\right){\rm d}s.\]
In e.g.\ \cite[Prop. 3.2]{dave} it was shown that the following weak convergence holds:
\[N^{\alpha/2}  {\boldsymbol I}^{(N)}(\cdot) \to {\boldsymbol B}(\cdot),\]
where ${\boldsymbol B}(\cdot)$ denotes a zero-mean $d$-dimensional Brownian motion with covariance matrix
$ {\rm diag}\{{\boldsymbol \pi}\} D+D^{\rm T} {\rm diag}\{{\boldsymbol \pi}\}$. The fact that this matrix is nonnegative 
definite has been proven in \cite[Prop. 3.2]{mmou}, and hence it allows a Cholesky decomposion.
We, in addition, obtain the weak convergence of $ {\boldsymbol H}^{(N)}(\cdot)$, 
with $H_i^{(N)}(t):=\lambda_iI_i^{(N)}(t)$,
to a zero-mean $d$-dimensional Brownian motion
with covariance matrix (and thus also allowing a Cholesky decomposition) 
\begin{equation}
\label{ddd}
V:=\Lambda\left({\rm diag}\{{\boldsymbol \pi}\} D+D^{\rm T} {\rm diag}\{{\boldsymbol \pi}\}\right)\Lambda.
\end{equation}

It also follows that this term behaves essentially proportional to $N^{1-\beta-\alpha/2};$ more specifically,
it converges to a constant if $\alpha\leq 1$, and vanishes otherwise.
\item[$\circ$] We now consider the second term.
We note that
\[\frac{\rm d}{{\rm d}t} \langle \KK_i^{(N)}\rangle_t = 
N {\lambda}_iZ^{(N)}_i(t) + 
\mu_i M_i^{(N)}(t),\]
and consequently 
\[\frac{1}{N} \langle \KK_i^{(N)}\rangle_t = 
 {\lambda}_i\int_0^t Z^{(N)}_i(s){\rm d}s +
\mu_i \int_0^t \frac{ M_i^{(N)}(s)}{N}{\rm d}s,\]
which we can prove, using standard arguments (such as the ergodic theorem), to converge to
\[w_i(t):=\lambda_i \pi_it +\mu_i \int_0^t  \varrho\hII_i(1-e^{-\mu_i s}){\rm d}s.\]
We thus find that ${\KK}_i^{(N)}(\cdot)/\sqrt{N}$ converges to an appropriately scaled one-dimensional Brownian motion.

By doing similar steps for $ \langle \KK_i^{(N)}+\KK_j^{(N)}\rangle_t$, for $i\not=j$, we find that the quadratic covariation 
between ${\KK}_i^{(N)}(\cdot)/\sqrt{N}$ and ${\KK}_j^{(N)}(\cdot)/\sqrt{N}$ equals 0.
We conclude from the above that ${\boldsymbol \kappa}^{(N)}(\cdot)/\sqrt{N}$ converges to an appropriately scaled
$d$-dimensional Brownian motion; the variance of component $i$ at time $t$ is $w_i(t)$, and the covariances are all 0. 

It also follows that this term is essentially proportional to $N^{1/2-\beta}$, 
which  is a constant if $\alpha\geq 1$, and vanishes otherwise.
\end{itemize}

Define $\tilde {\bs M}(t):=(\tilde M_1(t),\ldots, \tilde M_d(t))^{\rm T}$; we also write \[
W_i(t):= {\lambda_i\pi_i +\mu_i \varrho\hII_i(t)}= {2\lambda_i \pi_i-\lambda_i \pi_i e^{-\mu_i t}}.\]
Based on the above, we obtain that $\tilde M^{(N)}$ converges as $N\to\infty$ to the solution $\tilde {\bs M}(t)$ of the
stochastic differential equation, in self-evident notation,
\[{\rm d}\tilde {\bs M}(t) =
{-}{\mathcal M}\, \tilde {\bs M}(t) {\rm d}t+\,
{\sqrt{V1_{\{\alpha\le 1\}}+{\rm diag}\{{\bs W}(t)\}1_{\{\alpha \ge 1\}}}}\,{\rm d}{\boldsymbol B}(t).\]
From this stochastic differential equation it follows by applying standard techniques that the resulting limiting process
is a centered Gaussian process, with
\begin{eqnarray*}{\mathbb C}{\rm ov}\left(\tilde M_i(t),\tilde M_j(t)\right) &= &
e^{-\mu_i t - \mu_j t}\int_0^t e^{\mu_i s+\mu_j s}\lambda_i\lambda_j\left( \pi_i D_{ij} +\pi_{{j}} D_{ji}\right) {\rm d}s\\
&=&\frac{\lambda_i\lambda_j}{\mu_i+\mu_j}\left(1-e^{-(\mu_i+\mu_j)t}\right)\left( \pi_i D_{ij} +\pi_j D_{ji}\right)
\end{eqnarray*}
if $\alpha\le  1$. If $\alpha\ge  1$, on the contrary, the  covariance is $0$ if $i\not=j$, and $\varrho\hII_i(t)$ 
if $i=j.$ 

Now consider the limiting distribution of the {\it total} population of the system; from the above, it immediately 
follows that we have the weak convergence
\[\sum_{i=1}^d \tilde M^{(N)}_i(t)\to \tilde M(t):=\sum_{i=1}^d \tilde M_i(t),\]
which is a centered Gaussian process  of the Ornstein-Uhlenbeck type, characterized by its covariance $v\hII(t,u)$, 
as given in (\ref{limcII}). 
\begin{theorem} Consider Model {\sc ii}. As $N\to\infty$, the process $\tilde M^{(N)}(\cdot)$ converges weakly to a 
centered Gaussian process with covariance structure $v\hII(\cdot,\cdot)$ given
in \eqref{eqcov}.
\end{theorem}

\section{Numerical experiments}

In this section, we illustrate the results with two plots. In all cases, we consider a Model {\sc i} scenario with a two-state Markov chain. We assume that $q_{12}=q_{21}=5$, and ${\boldsymbol\lambda} =[20,\,10]$.

\begin{figure}
    \centering
    \includegraphics[width=0.7\linewidth]{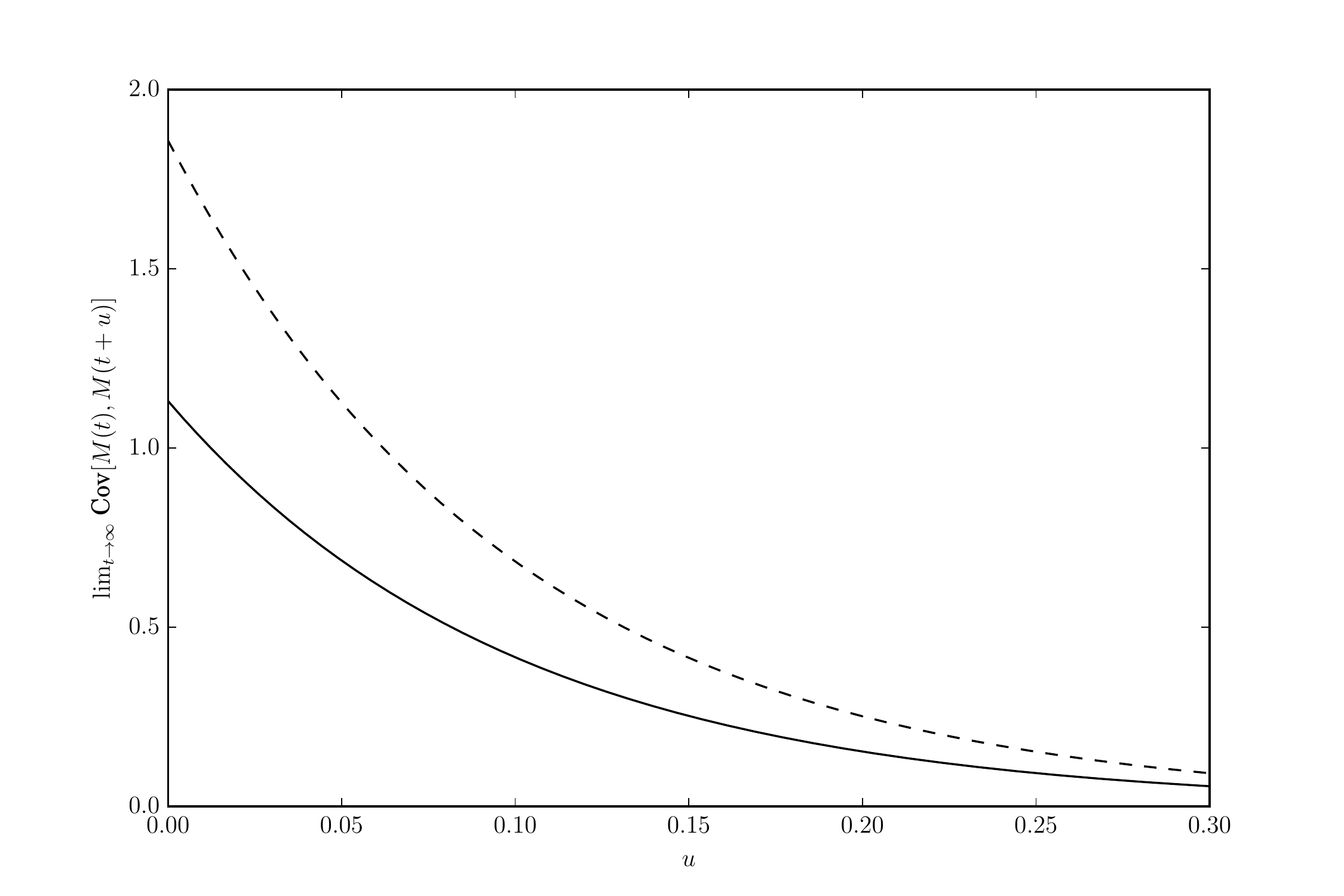}
    \caption{\label{figcov}{The limiting covariance versus time for ${\boldsymbol\mu} =[2,\, 1]$ (dashed line), and, ${\boldsymbol\mu} =[1,\,2]$ (solid line).}}
\end{figure}

In the first figure, Fig.~\ref{figcov}, we plot the covariance of a system starting in stationarity for two scenarios. In the first scenario (dashed line) ${\boldsymbol\mu} =[2,\,1]$, whereas in the second scenario (full line) ${\boldsymbol\mu} =[1, \,2]$.

\begin{figure}[t]
    \centering
    \includegraphics[width=0.7\linewidth]{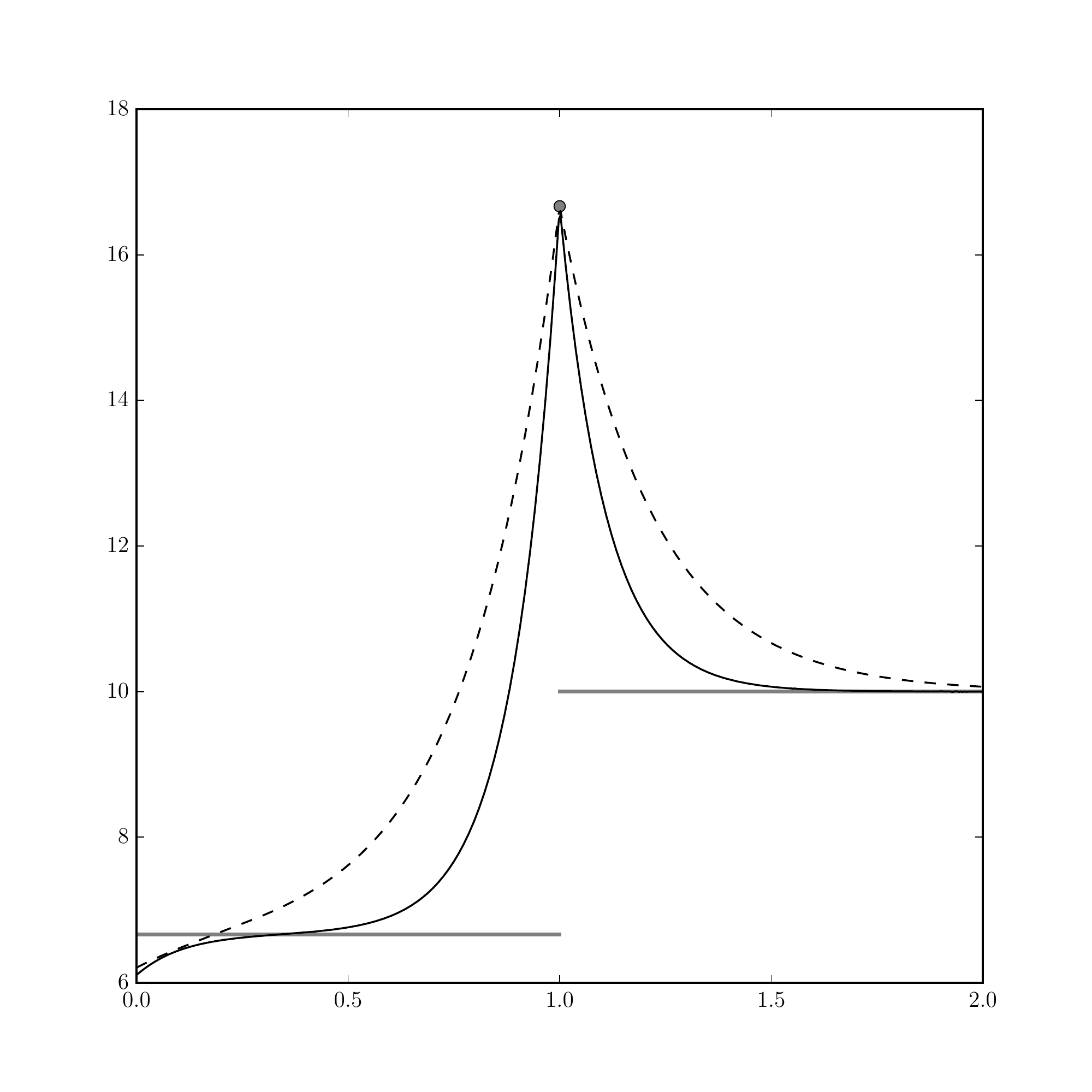}
    \caption{\label{figvar}{The scaled variance of $M^{(N)}$; dashed lines: $N=100$, solid: $N=100\,000$, limit $M(t)$: gray}}
\end{figure}

For the second plot, we assume ${\boldsymbol\mu} =[1\,,\, 2]$ and apply the scaling ${\boldsymbol\lambda} \mapsto N {\boldsymbol\lambda}$ and $Q\mapsto N^\alpha Q$.  Fig.~\ref{figvar} shows the stationary variance of the number of jobs. We divide this variance by the theoretically predicted growth factor $N^{2\gamma}$, and plot it against $\alpha$. The dashed line corresponds to $N=100$, the full line to $N=100\,000$.  We plot the limit $M(t)$ in gray.

\subsection*{Acknowledgment}
The authors thank Peter Spreij (Korteweg-de Vries Institute for Mathematics, University of Amsterdam) for useful suggestions.

\section*{Appendix A}
In this appendix we characterize the probability generating functions $\Xi_{ij}(\cdot,\cdot,\cdot,\cdot)$ by setting up a 
system of partial differential equations. 
The starting point for this is the system of Kolmogorov equations related to the transient probabilities  of the number 
of jobs present (jointly with the background state) at two time epochs:
\[p_{ij}(m,n,t,u):={\mathbb P}(M(t) = m, M(t+u)= n, J(t)=i, J(t+u)=j);\]
we suppress the $u$ as this is held fixed for the moment.
Standard arguments from Markov chain theory immediately yield the equations, 
for $i,j\in\{1,\ldots,d\}$, $m,n\in\{0,1,2\ldots\}$, {and $q_i:=-q_{ii}$},
\begin{eqnarray*}
p_{ij}(m,n,t+\Delta t,u)&=&p_{ij}(m,n,t,u)\left(1- (\lambda_i+\lambda_j +m\,\mu_i+n\,\mu_j+q_i+q_j)\Delta t\right)\\
&&+\,p_{ij}(m-1,n,t,u)\lambda_i \Delta t\:+\: p_{ij}(m,n-1,t,u)\lambda_j \Delta t\\
&&+\,p_{ij}(m+1,n,t,u)\,(m+1)\mu_i \Delta t\:+\: p_{ij}(m,n+1,t,u)\,(n+1)\mu_j \Delta t\\
&&+\,\sum_{k\not=i}p_{kj}(m,n,t,u)q_{ki} \Delta t\:+\:\sum_{k\not=j} p_{ik}(m,n,t,u)q_{kj} \Delta t +o(\Delta t);
\end{eqnarray*}
here $p_{ij}(-1,n,t,u)$ and $p_{ij}(m,-1,t,u)$ are to be understood as $0$.
As a consequence, with $p_{ij}'(m,n,t,u)$ denoting the derivative of $p_{ij}(m,n,t,u)$ with respect to $t$, it is readily 
obtained that
the transient probabilities satisfy the following system of (ordinary) differential equations:
\begin{eqnarray*}\lefteqn{p_{ij}'(m,n,t,u)=\lambda_i \left(p_{ij}(m-1,n,t,u)-p_{ij}(m,n,t,u)\right)+
\lambda_j\left(p_{ij}(m,n-1,t,u)-p_{ij}(m,n,t,u)\right)}\\
&&+\,\mu_i\left((m+1)p_{ij}(m+1,n,t,u)-mp_{ij}(m,n,t,u)\right)+
\mu_j\left((n+1)p_{ij}(m,n+1,t,u)-np_{ij}(m,n,t,u)\right)
\\
&&+\,\sum_{k=1}^dp_{kj}(m,n,t,u)q_{ki} \:+\:\sum_{k=1}^d p_{ik}(m,n,t,u)q_{kj} .
\end{eqnarray*}
Our goal is to transform these differential equations into a system of partial differential equations for the 
corresponding probability generating functions. 
To this end, multiply both sides of the equation by $z^mw^n$, and sum over $m,n=0,1,2,\ldots$. This results in the
following equation:
\begin{eqnarray*}
\frac{\partial  }{\partial t} \Xi_{ij}(z,w,t,u)&=&\left(\lambda_i(z-1)+\lambda_j(w-1) \right)\cdot \Xi_{ij}(z,w,t,u)\\
&&-\, \mu_i (z-1) \frac{\partial }{\partial z} \Xi_{ij}(z,w,t,u)- \mu_j(w-1) \frac{\partial }{\partial w} \Xi_{ij}(z,w,t,u)\\
&&+\,\sum_{k=1}^d\Xi_{kj}(z,w,t,u)q_{ki} \:+\:\sum_{k=1}^d \Xi_{ik}(z,w,t,u)q_{kj} ,
\end{eqnarray*}
which in matrix notation coincides with Eqn.\ (\ref{xi}).

\section*{Appendix B}
As will be proven in Appendix C below, the statement (\ref{lime}) can be refined to, for some constant $\kappa$ (whose precise form is 
irrelevant here),
\begin{equation}\label{appres}
{\mathbb E} M\hN(t) =( N\,\varrho\hI + N^{1-\alpha}\kappa)(1-e^{-\mu_\infty t})  +o(1).
\end{equation}
Likewise, for any $x\ge 0$,
\begin{equation}\label{asm}{\mathbb E} \left(M\hN(t)\,|\,M\hN(0)=x\right)
=(N\varrho\hI+N^{1-\alpha}\kappa)(1-e^{-\mu_\infty t})+
x e^{-\mu_\infty t}+o(1).\end{equation}
In the sequel we write $a(N):=N\varrho\hI+N^{1-\alpha}\kappa$.
Applying an elementary time shift, we obtain
that
 \begin{eqnarray*}
{\mathbb E}\left( M\hN(t)M\hN(t+u)\right) &=&
\int_0^\infty {\mathbb E} \left(M\hN(t)M\hN(t+u)\,|\,M\hN(t)=x\right){\mathbb P}(M\hN(t)\in{\rm d} x)\\
&=&\int_0^\infty x{\mathbb E} \left(M\hN(u)\,|\,M\hN(0)=x\right){\mathbb P}(M\hN(t)\in{\rm d}  x).
\end{eqnarray*}
By plugging in (\ref{asm}), the expression in the last display equals
\begin{eqnarray*}
\lefteqn{\int_0^\infty x\left(a(N)(1-e^{-\mu_\infty u})+
x e^{-\mu_\infty u}+o(N^{1-\alpha})\right){\mathbb P}(M\hN(t)\in{\rm d}  x)}\\
&=&\left(a(N)(1-e^{-\mu_\infty u})+o(N^{1-\alpha})\right){\mathbb E}M\hN(t)
+e^{-\mu_\infty u}\,{\mathbb E} \left((M\hN(t))^2\right)\\
&=&(\xi_N(u)+o(N^{1-\alpha}))(\xi_N(t)+o(N^{1-\alpha}))
+e^{-\mu_\infty u}\,{\mathbb E} \left((M\hN(t))^2\right),
\end{eqnarray*}
where $\xi_N(u):=a(N)(1-e^{-\mu_\infty u}).$
Now note that, due to the computations underlying \cite[Thm.~2]{PEIS},
with $\varrho\hI(t)$ and
$\varsigma\hI(t)$ as defined in Section \ref{SR},
and with $\psi(N)=o(N^{\max\{1,2-\alpha\}}),$ that
\[{\mathbb E} \left((M\hN(t))^2\right)
-\left({\mathbb E} M\hN(t) \right)^2= 
N^{2-\alpha}\varsigma\hI(t)1_{\{\alpha\le 1\}}+N\varrho\hI(t)1_{\{\alpha\ge 1\}}+\psi(N).\]
We now turn our attention to characterizing the covariance between $M\hN(t)$ and $M\hN(t+u)$. 
Based on the above  {we }have
\begin{eqnarray*}
\lefteqn{\hspace{-1.7cm}\COV(M\hN(t),M\hN(t+u)) =
{\mathbb E}\left(M\hN(t)M\hN(t+u)\right) - {\mathbb E}\left(M\hN(t)\right){\mathbb E}\left(M\hN(t+u)\right)}\\
&=&(\xi_N(u)+o(N^{1-\alpha}))(\xi_N(t)+o(N^{1-\alpha}))+e^{-\mu_\infty u}\left(\xi_N(t)+o(N^{1-\alpha})\right)^2\\ 
&&+\,e^{-\mu_\infty u}\left(N^{2-\alpha}\varsigma\hI(t)1_{\{\alpha\le 1\}} +
N\varrho\hI(t)1_{\{\alpha\ge 1\}}+\psi(N)\right)\\
&&-\,\left(\xi_N(t)+o(N^{1-\alpha})\right)\left(\xi_N (t+u))+o(N^{1-\alpha})\right).
\end{eqnarray*}
A direct computation now yields that {the zero-order} terms cancel, and that we end up with (\ref{limc}).

\section*{Appendix C}

In this appendix, we establish (\ref{appres}). 
The idea is  to manipulate  differential equation  (\ref{mteq}), so as to characterize the behavior
of $ {\bs m}(t) $ for $N$ large. 
The first step is to postmultiply the equation by the fundamental matrix $F:=D+\Pi$. We obtain, multiplying the equation
by $N^{-\alpha}$ as well,
\[{\bs m}(t) = {\bs m}(t) \Pi -N^{-\alpha}\,{\bs m}(t) {\mathcal M}F +
N^{1-\alpha}\,{\bs\pi}^{\rm T}\Lambda F-N^{-\alpha}{\bs m}'(t)F.\]
Iterate this equation once, we obtain
\begin{eqnarray*}\vm(t)&=&\vm(t)\Pi-N^{-\alpha}\vm(t){\mathcal M}\Pi+
N^{1-\alpha}{\bs \pi}^{\rm T}\Lambda\Pi-N^{-\alpha}\vm'(t)\Pi\\&&-\,N^{-\alpha}\vm(t)\Pi {\mathcal M}F-
N^{1-2\alpha}{\bs \pi}^{\rm T}\Lambda F{\mathcal M}F +N^{1-\alpha}{\bs\pi}^{\rm T}\Lambda F-N^{-\alpha}\vm'(t)\Pi
+o(N^{-\alpha}).\end{eqnarray*}
{Iterating once again to replace all occurrences of ${\bs m}(t)$ by ${\bs m}(t) \Pi$}, 
we obtain, with ${\bs n}(t):={\bs m}(t)\Pi$,
\begin{eqnarray*}\vm(t)&=&
{\bs n}(t)-N^{-\alpha}{\bs n}(t)\Pi {\mathcal M}\Pi\ -N^{1-2\alpha}{\bs \pi}^{\rm T}\Lambda F {\mathcal M}\Pi
+N^{1-\alpha}{\bs \pi}^{\rm T}\Lambda\Pi-N^{-\alpha}{\bs n}'(t)\\
&&-\,N^{-\alpha}{\bs n}(t) \Pi {\mathcal M}F-
N^{1-2\alpha}{\bs \pi}^{\rm T}\Lambda F{\mathcal M}F +N^{1-\alpha}{\bs\pi}^{\rm T} \Lambda F-N^{-\alpha}{\bs n}'(t)+
o(N^{-\alpha}).\end{eqnarray*}
Now postmultiply this equation by $N^\alpha\Pi {{\bs 1}}$.
Recalling that $F{\bs 1}={\bs 1}$ and $\Pi{\bs 1}={\bs 1}$, 
we obtain
\[{\bs n}'(t){\bs 1}=-{\bs n}(t)\Pi {\mathcal M}{\bs 1}+N{\bs\pi}^{\rm T}\Lambda{\bs 1}
 - N^{1-\alpha}{\bs \pi}^{\rm T} \Lambda F{\mathcal M}{\bs 1}+o(1).\]
which leads to, using $ {\bs n}(t){\bs 1}:=\phi(t)$ and $\Pi={\bs 1}{\bs\pi}^{\rm T}$,
\[\phi'(t)=-\mu_\infty\phi(t)+N\lambda_\infty  
- N^{1-\alpha}{\bs \pi}^{\rm T} \Lambda F{\mathcal M}{\bs 1}+o(1).\]
We find that, with $\phi(0)=0$,
\[{\mathbb E}M\hN(t) = \left(N\frac{\lambda_\infty}{\mu_\infty}  
-\frac{1}{\mu_\infty} N^{1-\alpha}{\bs \pi}^{\rm T} \Lambda F{\mathcal M}{\bs 1}\right)
(1-e^{-\mu_\infty t})+o(1).\]

\section*{Appendix D}
We first focus on the first term in the right hand side of (\ref{TC}). To this end, consider the following decomposition:
\[M(t):= M^{(1)}(t,t+u)+ M^{(2)}(t,t+u),\:\: \:\:M(t+u):= M^{(2)}(t,t+u)+ M^{(3)}(t,t+u),\]
where $M^{(1)}(t,t+u)$ are the jobs that arrived in $[0,t)$ that are still present at time $t$ 
but have left at time $t+u$, $M^{(2)}(t,t+u)$ the jobs that have arrived in 
$[0,t)$ that are still present at time $t+u$, and $M^{(3)}(t,t+u)$ the jobs that have arrived in $[t,t+u)$ 
that are still present at time $t+u$.
 Observe that, conditional on $J$,
these three random quantities are independent. As a result,
\[\MEAN(\COV(M(t),M(t+u))\,|\,J)) = \MEAN(\VAR \,M^{(2)}(t,t+u)\,|\,J).\]
Mimicking the arguments used in \cite{DAURIA}, it is immediate that 
$M^{(2)}(t,t+u)$, conditional on $J$, has a Poisson distribution with parameter
\[\int_0^t \lambda_{J(s)} e^{-\mu_{J(s)}(t+u-s)}{\rm d}s.\]
We conclude that
\begin{eqnarray*}
\MEAN(\COV(M(t),M(t+u))\,|\,J))&=&\MEAN\left(\int_0^t \lambda_{J(s)} e^{-\mu_{J(s)}(t+u-s)}{\rm d}s \right) \\&=& 
\sum_{i=1}^d \pi_i \lambda_i\int_0^t e^{-\mu_i(t+u-s)} {\rm d}s
=
\sum_{i=1}^d \pi_i \frac{\lambda_i}{\mu_i}\left(1-e^{-\mu_i t}\right)e^{-\mu_iu}
.\end{eqnarray*}
Now analyze the second term in the right hand side of (\ref{TC}). First observe that it can be written as 
\[\COV\left(\int_0^t \lambda_{J( r)} e^{-\mu_{J( r)}(t-r)}{\rm d}r, 
\int_0^{t+u}\lambda_{J(s)} e^{-\mu_{J(s)}(t+u-s)}{\rm d}s  \right).\]
This decomposes into 
$I_1+I_2$, where
\begin{eqnarray*}
I_1&:=&\sum_{i=1}^d\sum_{j=1}^d\lambda_i\lambda_j\cK_{ij},\:\:\mbox{where}\:\:\cK_{ij}:=
\int_0^t\int_0^s e^{-\mu_i(t-r)} e^{-\mu_j(t+u-s) }\pi_i\left(p_{ij}(s-r)-\pi_j\right){\rm d}r{\rm d}s,\\
I_2&:=&\sum_{i=1}^d\sum_{j=1}^d\lambda_i\lambda_j\cL_{ij},\:\:\mbox{where}\:\:\cL_{ij}:=
\int_0^t\int_s^{t+u} e^{-\mu_i(t-r)} e^{-\mu_j(t+u-s) }\pi_j\left(p_{ji}(r-s)-\pi_i\right){\rm d}r{\rm d}s.
\end{eqnarray*}
Let us first evaluate $\cK_{ij}\equiv \cK_{ij}(t,u).$ To this end, substitute $w:=s-r$ (i.e., replace $r$ by $s-w$), 
and then 
interchange the order of integration, so as to obtain
\[\cK_{ij}= e^{-\mu_j(t+u) }\pi_i
\int_0^t\left(\int_w^t e^{(\mu_i+\mu_j)s}{\rm d}s\right) e^{-\mu_i(t+w)} \left(p_{ij}(w)-\pi_j\right){\rm d}w.\]
Performing the inner integral (i.e., the one over $s$) leads to
\[\cK_{ij} =\frac{1}{\mu_i+\mu_j} 
e^{-\mu_j(t+u)} \pi_i \int_0^t \left(e^{-\mu_iw+\mu_jt}-e^{-\mu_it+\mu_jw}\right)\left(p_{ij}(w)-\pi_j\right){\rm d}w.\]
For $\cL_{ij}\equiv \cL_{ij}(t,u)$, again by a substitution and by interchanging the order of integration, we obtain
$\cL_{ij}^{(1)}+\cL_{ij}^{(2)}$, where
\begin{eqnarray*}
\cL_{ij}^{(1)}&:=&e^{-\mu_j(t+u)}\pi_j \int_0^{u}\left(\int_0^t e^{(\mu_i+\mu_j)s}{\rm d}s\right)
e^{-\mu_i(t-w)} \left(p_{ji}(w)-\pi_i\right){\rm d}w
,\\
\cL_{ij}^{(2)}&:=&e^{-\mu_j(t+u)}\pi_j \int_{u}^{t+u}\left(\int_0^{t+u-w}e^{(\mu_i+\mu_j)s}{\rm d}s\right)
e^{-\mu_i(t-w)} \left(p_{ji}(w)-\pi_i\right){\rm d}w,
\end{eqnarray*}
which reduce to 
\begin{eqnarray*}
\cL_{ij}^{(1)}&:=&\frac{1}{\mu_i+\mu_j}e^{-\mu_j(t+u)}\pi_j 
\left(e^{\mu_j t}-e^{-\mu_i t}\right)
\int_0^{u} e^{\mu_i w}
 \left(p_{ji}(w)-\pi_i\right){\rm d}w
,\\
\cL_{ij}^{(2)}&:=&\frac{1}{\mu_i+\mu_j}e^{-\mu_i t}\pi_j 
\int_{u}^{t+u}
\left(e^{\mu_i (t+u) -\mu_jw}- e^{\mu_i w-\mu_j(t+u)}\right) \left(p_{ji}(w)-\pi_i\right){\rm d}w.
\end{eqnarray*}
Now Eqn.\ (\ref{eqcov}) follows.

\end{document}